\documentclass[12pt]{amsart}
\usepackage{amsmath}
\usepackage{amssymb}
\def\r{\mathrm}

\def\e{\varepsilon}

\def\Z{\mathbb Z}
\def\R{\mathbb R}

\def\eqdef{{\buildrel \r{def} \over =}}

\def\oo{\infty}

\def\half{{\tfrac12}}

\def\upref#1{\upn{\ref{#1}}}

\theoremstyle{definition}

\theoremstyle{plain}
\newtheorem{thm}{Theorem}
\newtheorem{cor}{Corollary}
\newtheorem{lem}{Lemma}
\newtheorem{prp}{Proposition}
\theoremstyle{remark}
\newtheorem{rmrk}{Remark}
\newtheorem{ex}{Example}

\title[A general Hsu-Robbins-Erd\"os type estimate]{A general
Hsu-Robbins-Erd\"os type estimate of tail probabilities
of sums of independent identically distributed random variables\\
{\em Une estim\'ee g\'en\'erale de type Hsu-Robbins-Erd\"os pour les probabilit\'es 
des queues des sommes des variables al\'eatoires ind\'ependantes
et identiquement distribu\'ees}}
\author{Alexander R. Pruss}
\date{October 18, 1998}
\address{Department of Philosophy\\
University of Pittsburgh\\
Pittsburgh, PA 15260\\
U.S.A.}
\email{pruss+@pitt.edu}
\begin{document}
\begin{abstract}
        Let $X_1,X_2,\dots$ be a sequence of independent and identically
        distributed random variables, and put $S_n=X_1+\dots+X_n$.
        Under some conditions on the positive sequence $\tau_n$ and the
        positive increasing sequence $a_n$, we
        give necessary and sufficient conditions for the convergence of
$\sum_{n=1}^\oo \tau_n P(|S_n|\ge \e a_n)$
        for all $\e>0$, generalizing Baum and Katz's~(1965)
        generalization of the Hsu-Robbins-Erd\H os (1947, 1949) law of
        large numbers, also allowing us to characterize the
        convergence of the above series in the case where
        $\tau_n=n^{-1}$ and $a_n=(n\log n)^{1/2}$ for $n\ge 2$, thereby
        answering a question of Sp\u ataru.  Moreover, some results for
        non-identically distributed independent random variables are
        obtained by a recent comparison inequality. Our basic method is to use a
        central limit theorem estimate of Nagaev~(1965) combined with
        the Hoffman-J\o rgensen inequality~(1974).

\medskip
{}
\noindent
{\sc R\'esum\'e.}.\ 
        Soit $X_1,X_2,\dots$ -- une suite des variables al\'eatoires
        ind\'ependantes et identiquement distribu\'ees, et met
        $S_n=X_1+\dots+X_n$. Dans certaines conditions sur la suite
        positive $\tau_n$ et la suite croissante positive $a_n$, nous
        donnons des conditions n\'ecessaires et suffisantes pour la
        convergence de $\sum_{n=1}^\oo \tau_n P(|S_n|\ge \e a_n)$ pour
        tout $\e>0$, g\'en\'eralisant l'extension de Baum et Katz~(1965)
        de la loi de grands nombres de Hsu-Robbins-Erd\H os (1947,
        1949).  Ce nous permet de caract\'eriser la convergence de la
        s\'erie ci-dessus avec $\tau_n=n^{-1}$ et $a_n=(n\log n)^{1/2}$
        pour $n\ge 2$, r\'epondant \'a une question de Sp\u ataru.
        D'ailleurs, nous obtenons quelques r\'esultats pour des 
        variables al\'eatoires ind\'ependantes mais non identiquement
        distribu\'ees, par une in\'egalit\'e r\'ecente de comparaison.
        Notre m\'ethode est bas\'ee sur une estim\'ee de Nagaev~(1965) dans
        le th\'eor\`eme de limite centrale, en combinaison avec
        l'in\'egalit\'e de Hoffman-J\o rgensen~(1974).
\end{abstract}

%%%%%%%%%%%%%%%%%%%%%%%%%%%%%%%%
\maketitle
\bibliographystyle{amsplain}
%%%%%%%%%%%%%%%%%%%%%%%%%%%%%%%%
\section{Introduction and main result}
Hsu and Robbins~\cite{HsuRobbins} have shown that if
$X_1,X_2,\dotsc$ are independent and identically distributed mean zero random variables with
finite variance, then \begin{equation}\label{eq:HR-concl}
        \sum_{n=1}^\infty P(|S_n| \ge \e n) < \infty, \qquad
        \forall \e >0,
 \end{equation}
where $S_n=X_1+\dots+X_n$.
Erd\H{o}s~\cite{Erdos,Erdos:Two} has shown
that a converse implication also holds, namely that if
\eqref{eq:HR-concl} holds  and the $\{ X_i \}$ are
independent and identically distributed, then $E[X_1^2]<\infty$ and $E[X_1]=0$.
Since then, a number of extensions in several directions have been
proved; see \cite{LRJW} and \cite{Pruss:Spataru} for partial
bibliographies and brief discussions.  The purpose of the present paper
is to prove a new rather general extension of the Hsu-Robbins-Erd\H os
law of large numbers.

Throughout, terms like ``positive'' and ``increasing''
indicate the non-strict varieties (``non-negative'' and
``non-decreasing'', respectively).

We say that a positive sequence $\{ \tau_n \}_{n=1}^\oo$ {\em satisfies 
Condition~A} provided that for all positive decreasing
series $\{ c_n \}_{n=1}^\oo$
such that $\sum_{n=1}^\oo \tau_n
\min(nc_n,1)$ converges, we likewise have $\sum_{n=1}^\oo \tau_n nc_n$
converging.  The following Proposition gives us a simple sufficient criterion
for Condition~A.  A proof will be given in
Section~\ref{sec:proofs}, below.

\begin{prp}\label{prp:condA} Suppose that either
$\liminf_{n\to\oo}\tau_n>0$ or that both of the following two conditions
are satisfied:
\begin{enumerate}
\item[\upn{(a)}] $\liminf_{n\to\oo}n\tau_n > 0$, and
\item[\upn{(b)}] there is a constant $C\in (0,\oo)$ such that for each
$n\in\Z^+_0$, if\/ $2^{n-1}\le k<2^n$, then
$$
        C\tau_{2^{n-1}} \ge \tau_k \ge C^{-1} \tau_{2^n}.
$$
\end{enumerate}
Then, $\{ \tau_n \}$ satisfies Condition~\upn{A}.
\end{prp}

We now state our main result.
\begin{thm}\label{th:main} Let $X_1,X_2,\dotsc$ be independent
identically distributed random variables.
Put $S_n=X_1+\dots+X_n$.
Let $\{ \tau_n \}_{n=1}^\oo$ be any positive sequence of numbers,
and let $\{ a_n \}_{n=1}^\oo$ be an increasing strictly positive
sequence tending to $\oo$.
Furthermore, suppose there exist finite real constants
$\theta\ge 1$, $C>0$ and $N\ge 2$ such that for all $n\ge N$ we have:
\begin{equation}\label{eq:aux-cond}
        \frac{a_n^{3\theta}}{n^{\theta-1}} \sum_{k=n}^\oo
        \frac{k^{\theta}\tau_k}{a_k^{3\theta}}
          \le C \sum_{k=1}^{n-1} k\tau_k.
\end{equation}
If we have $\theta>1$, then additionally assume that
\begin{equation}\label{eq:aux-aux-cond}
        \liminf_{n\to\oo} \inf_{k\ge n} \frac{a_k^3}{k a_n^3}
           \sum_{j=1}^{n-1} j\tau_j > 0.
\end{equation}

Then, if all of the following three conditions hold:
\begin{enumerate}
\item[\upn{(i)}]
        there is a sequence $\{ \mu_n \}_{n=1}^\oo$ with $\mu_n$ a
        median of $S_n$ for all $n$, and such that for all $\e>0$,
        we have $\sum_{n\in M}
        \tau_n<\oo$, where $M=\{ n\in \Z^+ : |\mu_n| > \e a_n \}$,
\item[\upn{(ii)}] $\sum_{n=1}^\oo n\tau_n P(|X_1|\ge \e a_n)<\oo$ for
all $\e>0$, and
\item[\upn{(iii)}] $\sum_{n=1}^\oo \tau_n
e^{-\e^2 a_n^2 / (nT_{\e,n})} <\oo$ for all $\e>0$,
where $T_{\e,n}=E[X_1^2\cdot 1_{\{|X_1|<\e a_n\}}]$,
and where $e^{-t/0}=0$ for all $t>0$,
\end{enumerate}
we will have:
\begin{equation}\label{eq:conv}
        \sum_{n=1}^\oo \tau_n P(|S_n|\ge \e a_n) < \oo, \quad \forall
        \e>0.
\end{equation}
Conversely, if the sequence $\{ \tau_n \}$ satisfies Condition~A and
\eqref{eq:conv} holds, conditions~\upn{(i)}, \upn{(ii)} and~\upn{(iii)}
will also hold.
\end{thm}

The proof will be given in Section~\ref{sec:proofs}, below.  The result
is closely related to work of Klesov~\cite[Theorem~4]{Klesov}, though we
are working with a more general class of sequences $\{ a_n \}$.  Our
proof will be based upon Klesov's Hoffman-J\o rgensen inequality based
approach, combined with a central limit theorem estimate of
Nagaev~\cite{Nagaev}, the latter being used rather like in
\cite{Pruss:Riemann}.

The particular newness of the result is that it works for $a_n$ near
the critical growth $n^{1/2}$ involved in the central limit theorem.
For instance in the next section (see Corollary~\ref{cor:Sp}) we will
use
the Theorem to characterize the convergence of \eqref{eq:conv} in the
case where $\tau_n=n^{-1}$ and $a_n = (n\log n)^{1/2}$ (for $n\ge 2$),
thereby answering a question of Sp\u ataru.

\begin{rmrk}\label{rk:sv-1}
Recall that a measurable function $\phi$ on $[0,\oo)$ is slowly varying (in
the sense of Karamata) providing that for all $\lambda>0$ we have
$\lim_{x\to\oo} \phi(\lambda x)/\phi(x)=1$.
It is not difficult to see that if $L$ and $K$ are strictly positive
slowly varying functions on $[0,\oo)$, and if $a_n=L(n) n^\alpha$ for some
$\alpha>\tfrac13$ while $\tau_n=K(n) n^\beta$ for some $\beta$ satisfies
$\liminf n\tau_n > 0$, then
conditions \eqref{eq:aux-cond} and \eqref{eq:aux-aux-cond} will be satisfied
for a sufficiently large $\theta$.  To see this, use the condition that
$\liminf k\tau_k > 0$ to observe that $n\le c\sum_{k=1}^{n-1} k\tau_k$
for all $n\ge N$, for some constants $c$ and $N$, with $c$ independent
of $n$.  Then, for sufficiently
large $\theta$, approximate $\sum_{k=n}^\oo a_k^{-3\theta} k^\theta
\tau_k$ by an integral (it is here that the condition that
$\alpha>\tfrac13$ will be used;  one will also need to use the easy fact
that powers and products of slowly varying functions are slowly varying)
and use the fact that if $\phi$ is slowly varying then $\int_{x}^\oo
t^{-p} \phi(t)\, dt$ is asymptotic to $(p-1)^{-1} x^{1-p} \phi(x)$ for
large $x$ by \cite[Proposition~1.5.10]{BGT}. This will yield
\eqref{eq:aux-cond}.  Condition~\eqref{eq:aux-aux-cond} can be obtained
by noting that $L$ can without loss of generality be replaced by a
normalized slowly varying function (see \cite[p.~15]{BGT}) asymptotic to
$L$, and then $a_k^3/k$ will be an increasing function for sufficiently
large $k$ if $\alpha> 1/3$, by the Bojanic-Karamata
theorem~\cite[Theorem~1.5.5]{BGT}, so that the infimum on the left hand
side of \eqref{eq:aux-aux-cond} will be attained at $k=n$, and the
truth of \eqref{eq:aux-aux-cond} will follow from the condition that
$\liminf n\tau_n > 0$.  Condition~A will also be satisfied by $\{ \tau_n
\}$ under the above circumstances, as is very easy to see by Proposition~\ref{prp:condA}.
\end{rmrk}

\begin{rmrk}\label{rk:MS}  By a maximal inequality of
Montgomery-Smith~\cite{MS}, condition~\eqref{eq:conv} is equivalent to: $$
\sum_{n=1}^\oo \tau_n P\bigl(\bigl|\sup_{1\le k\le n} S_k\bigr|\ge 
\e a_n\bigr)<\oo,
\qquad\forall{\e>0}. $$
\end{rmrk}

\begin{rmrk}
Readers familiar with Hsu-Robbins-Erd\H os laws of large numbers may be
surprised at condition~(iii) in Theorem~\ref{th:main}, since normally
these laws of large numbers simply have \eqref{eq:conv} equivalent to
(i) and (ii) (under appropriate conditions on $\{ \tau_n \}$ and $\{ a_n \}$).
In general, condition~(iii) cannot be eliminated from
Theorem~\ref{th:main}, and seems to become particularly significant for
$a_n$ near the critical grown $n^{1/2}$.
(Example~\ref{ex:MS} in Section~\ref{sec:cor} will
show that condition~(iii) cannot be eliminated if $\tau_n=n^{-1}$ and
$a_n=(n\log n)^{1/2}$.)  However, in a number of special cases, 
condition~(iii) can indeed be removed, as is seen in the following
result.
\end{rmrk}

\begin{thm}\label{th:aux}  Suppose the $\{ \tau_n \}$ is positive and $\{ a_n \}$ is increasing
and strictly positive, and that there are finite real constants
$\theta\ge 1$, $C>0$ and $N\ge 2$ such that for all $n\ge N$ we have:
\begin{equation}\label{eq:aux-cond2}
        \frac{a_n^{2\theta}}{n^{\theta-1}} \sum_{k=n}^\oo
        \frac{k^{\theta}\tau_k}{a_k^{2\theta}}
          \le C \sum_{k=1}^{n-1} k\tau_k
\end{equation}
and
\begin{equation}\label{eq:aux-aux-cond2}
        \liminf_{n\to\oo} \inf_{k\ge n} \frac{a_k^2}{k a_n^2}
           \sum_{j=1}^{n-1} j\tau_j > 0.
\end{equation}
Then, any random variable $X_1$ satisfying condition~\upn{(ii)} of
Theorem~\upref{th:main} automatically satisfies condition~\upn{(iii)} of
that Theorem.
\end{thm}

The proof will be given in Section~\ref{sec:proofs}, below.  Note that
\eqref{eq:aux-cond2} automatically implies \eqref{eq:aux-cond}, and
\eqref{eq:aux-aux-cond2} likewise implies \eqref{eq:aux-aux-cond}.

\begin{rmrk}\label{rk:sv-2}
It is not difficult to see that if $L$ and $K$ are slowly varying
strictly positive functions on $[0,\oo)$, and if $a_n=L(n) n^\alpha$ for
some $\alpha> \tfrac12$ while $\tau_n=K(n) n^\beta$ for some
$\beta\in\R$ satisfies $\liminf n\tau_n>0$, then conditions \eqref{eq:aux-cond2} and
\eqref{eq:aux-aux-cond2} will be satisfied for a sufficiently large
$\theta$.  Since by Remark~\ref{rk:sv-1}, conditions~\eqref{eq:aux-cond} and
\eqref{eq:aux-aux-cond} would allow $\alpha$ to be $\half$ (or in fact any
value greater than $\tfrac13$), this helps to further illustrate how
condition~(iii) becomes relevant close to the critical growth $n^{1/2}$
of $a_n$,
but for faster growths (i.e., around $n^\alpha$ for $\alpha>\half$), it
can be eliminated by Theorem~\ref{th:aux}. \end{rmrk} 

\section{Corollaries and applications}\label{sec:cor}
We now obtain the following generalization due to
Baum and Katz~\cite{BaumKatz} of the Hsu-Robbins-Erd\H os
result, thereby showing that
Theorem~\ref{th:main} is indeed more general than
the Hsu-Robbins-Erd\H os law of large numbers. (See~\cite{LRJW} for a 
discussion of the pedigree of the Baum and
Katz result.)
\begin{cor}\label{cor:bk}
Suppose that $X_1,X_2,\dotsc$ are independent identically
distributed random variables.  Let $S_n=X_1+\dots+X_n$.
Fix $r\ge 1$ and $0<p<2$.
Then the conjunction of the conditions
\begin{enumerate}
\item[\upn{(a)}] if $p\ge 1$ then $E[X_1]=0$, and
\item[\upn{(b)}] $E[|X_1|^{rp}]<\oo$
\end{enumerate}
holds if and only if
\begin{equation}\label{eq:bk-concl}
        \sum_{n=1}^\oo n^{r-2} P(|S_n|\ge \e n^{1/p})<\oo, \qquad
        \forall{\e>0}.
\end{equation}
\end{cor}

\begin{proof}[Proof]
Let $\tau_n=n^{r-2}$ and $a_n=n^{1/p}$.
Observe that \eqref{eq:aux-cond}, \eqref{eq:aux-aux-cond},
\eqref{eq:aux-cond2}, \eqref{eq:aux-aux-cond2} and Condition~A all
hold for appropriate choices of
$\theta$ (one can use Remarks~\ref{rk:sv-1} and~\ref{rk:sv-2} here if
one so desires, but in fact a direct verification is easy).
Then by the Marcinkiewicz-Zygmund strong law of large numbers~\cite[\S 16.4.A.3]{Loeve},
if (a) and (b)
hold, then $S_n/a_n\to 0$ almost surely, hence also in probability, and
therefore condition~(i) of
Theorem~\ref{th:main} holds as well.  Moreover, (ii) is equivalent to (b),
and (ii) implies (iii) by Theorem~\ref{th:aux}.
Thus, by Theorem~\ref{th:main}, if (a)
and (b) hold, \eqref{eq:bk-concl} follows.  Conversely, if
\eqref{eq:bk-concl} holds, then since $\tau_n$ satisfies
Condition~A by Proposition~\ref{prp:condA}, it follows from
Theorem~\ref{th:main} that (i) and (ii)
(and (iii), but that is not needed) hold.  Hence (b) holds, since it is
equivalent to (ii).  It remains to show that (a) holds.
The easiest way to do this is to note that if $p\ge 1$ then by
the Marcinkiewicz-Zygmund strong law of large numbers, we have
$(S_n-nE[X_1])/a_n\to 0$ almost surely, hence in probability, and
therefore $(\mu_n-nE[X_1])/n^{1/p}\to 0$, where
$\mu_n$ is a median of $S_n$, and by condition~(i) it will then follow
that $E[X_1]=0$.
\end{proof}

Recall that random variables $X_1,\dots,X_n$ are said to be $K$-{\em
weakly
mean dominated} by a random
variable $X$ providing that for all $\lambda$ we have:
$$
        \frac1n \sum_{k=1}^n P(|X_k|\ge \lambda) \le K P(|X|\ge\lambda)
$$
(see \cite{Gut}).

\begin{cor}\label{cor:wmd}
Let $X$ be a random variable and $K$ any finite constant.
Let $\{ X_{kn} \}_{1\le k\le n}$ be a triangular array of
random variables such that
$X_{1n},\dots,X_{nn}$ are $K$-weakly mean dominated
by $X$ for each fixed $n$.
Let $\tau_n$ be a positive sequence of
numbers and let $a_n$ be an increasing strictly positive sequence for
$n\ge 1$.
Put $S_n=X_{1n}+\dots+X_{nn}$.
Suppose that \eqref{eq:aux-cond} and condition~\upn{(i)} of
Theorem~\upref{th:main} hold, and that
conditions~\upn{(ii)} and~\upn{(iii)} of
Theorem~\upref{th:main} are satisfied with $X$ in place of $X_1$.  Then,
$$
        \sum_{n=1}^\oo \tau_n P(|S_n|\ge \e a_n) < \oo
$$
for all $\e>0$.
\end{cor}

\begin{proof}
Without loss of generality we may assume that $X$ is symmetric (i.e.,\
$X$ and $-X$ have the same distribution), since otherwise we can replace
$X$ with $\alpha X$, where $\alpha$ is independent of $X$ and
$P(\alpha=1)=P(\alpha=-1)=\half$, and this replacement will not affect
the weak mean domination conditions.
By Theorem~\ref{th:main} we then have
\begin{equation}\label{eq:primed-conv}
        \sum_{n=1}^\oo \tau_n P(|S_n'|\ge \e a_n)<\oo,
        \qquad\forall\e>0,
\end{equation}
where $S_n'$ is the sum of $n$ independent copies of $X$ (condition~(i)
in this case will hold trivially by symmetry, while (ii) and (iii) were
assumed in the statement of Corollary~\upref{cor:wmd}.) One may slightly
modify the comparison result in \cite[Corollary~1]{Pruss:noniid} by
assuming our Theorem~\ref{th:main}'s condition~(i) in place of the 
assumption in that paper that $S_n/a_n\to 0$ in
probability, which modification only very slightly affects the proof
(one will need to use \eqref{eq:for-ref}, below, after obtaining the
convergence of \cite[series~(1.2)]{Pruss:noniid} in the original proof
in \cite{Pruss:noniid}).
Thus modified, \cite[Corollary~1]{Pruss:noniid} together with \eqref{eq:primed-conv}
yields the conclusion of our Corollary~\ref{cor:wmd}.
\end{proof}

The following Corollary yielding a result similar to one of Hu, Moricz and
Taylor~\cite{HMT} (cf.\ \cite{Gut,HSV}) can be derived from
Corollary~\ref{cor:wmd} exactly in the way that Corollary~\ref{cor:bk}
was derived from Theorem~\ref{th:main}.
\begin{cor}
Let $X$ be a random variable and let $K$ be any finite constant.
Let $\{ X_{kn} \}_{1\le k\le n}$ be a triangular array of
random variables such that $X_{1n},\dots,X_{nn}$ are $K$-weakly mean dominated
by $X$.  Let $S_n=X_{1n}+\dots+X_{nn}$.
Fix $r\ge 1$ and $0<p<2$.
Suppose $E[|X|^{rp}]<\oo$.  If $p\ge 1$ then assume also that
$E[X_{n1}+\dots+X_{nn}]=0$ for all $n$.
Then,
$$
        \sum_{n=1}^\oo n^{r-2} P(|S_n|\ge \e n^{1/p}), \qquad
        \forall{\e>0}.
$$
\end{cor}

Now define $\log^+ x=\log (2+x)$.  It is easy to see that
Theorem~\ref{th:main} implies the following result.
\begin{cor}\label{cor:Sp}
Let $X_1,X_2,\dotsc$ be independent and identically distributed random
variables.  Then,
\begin{equation}\label{eq:sp-conv}
        \sum_{n=2}^\oo n^{-1} P(|S_n| \ge \e (n\log n)^{1/2}) < \oo,
        \quad \forall \e>0,
\end{equation}
if and only if all of the following three conditions hold:
\begin{enumerate}
\item[\upn{(a)}] $E[X_1]=0$,
\item[\upn{(b)}] $E[X_1^2/\log^+ |X_1|]<\oo$, and
\item[\upn{(c)}] $\sum_{n=2}^\oo n^{-1-\e^2/T_{\e,n}}<\oo$ for
all $\e>0$, where $T_{\e,n}=E[X_1^2\cdot
1_{\{|X_1|<\e (n\log n)^{1/2}\}}]$.
\end{enumerate}
\end{cor}

\begin{rmrk} As per Remark~\ref{rk:MS}, we can replace $S_n$ by
$\max_{1\le k\le n} S_k$ in \eqref{eq:sp-conv}.
\end{rmrk}

The proof of the following Lemma will be given at the end of
Section~\ref{sec:proofs}.
\begin{lem}\label{lem:Sp}
        If condition \upn{(b)} of Corollary~\upn{\ref{cor:Sp}}
        holds, then $(S_n-E[S_n]) / (n\log n)^{1/2} \to 0$ almost surely, and
        hence also in probability, as $n\to\oo$.
\end{lem}

\begin{proof}[Proof of Corollary~\upref{cor:Sp}]
Let $\tau_n=1/n$ and $a_n=(n\log n)^{1/2}$ for $n\ge 2$.
Note that \eqref{eq:aux-cond} is easily seen to be
satisfied with $\theta=1$, and Condition~A holds by 
Proposition~\ref{prp:condA}.  Observe that (b) is
equivalent to condition~(ii) of Theorem~\ref{th:main} in the present
setting, and that if (b) holds,
then (a) is equivalent to (i) by
Lemma~\ref{lem:Sp}, so that the conjunction of conditions (a) and (b)
of Corollary~\ref{cor:Sp}
is
equivalent to that of conditions (i) and (ii) of Theorem~\ref{th:main}.  
Also, it is easy to check that (c) is
equivalent to (iii), since the summands in the sums in both conditions
are equal.
\end{proof}

Professor Aurel Sp\u ataru has asked the author whether
\eqref{eq:sp-conv} is equivalent to the conjunction of (a) and (b). This
would be expected by analogy with other Hsu-Robbins-Erd\H os laws of
large numbers (such as Corollary~\ref{cor:bk}). It
is this question that has inspired the present paper. In light of the
Corollary~\ref{cor:Sp}, Sp\u ataru's question is equivalent to asking
whether the conjunction of (a) and (b) implies (c).  The following
counterexample that Professor Stephen Montgomery-Smith has privately
communicated to the author  shows that the answer is negative, and hence so is
the answer to Sp\u ataru's question.

\begin{ex}\label{ex:MS}
Let $\{ K_m \}_{m=0}^\oo$ be a very rapidly increasing strictly positive
sequence, with $K_0=0$.
The degree of rapidity of increase will be chosen later so as to
make the argument go. Let $\psi(t)=(t \log t)^{1/2}$ for $t\ge 2$.
Extend $\psi$ linearly to the interval $[0,2]$ in such a way that $\psi(0)=0$.
Let $\phi$
be the inverse function of $\psi$.  Assume that $K_1\ge 2$. Let
$X_1$ be a random variable such that
$P(X_1=\psi(K_m))=P(X_1=-\psi(K_m))=2^{-m-1}/K_m$ for all $m$, and with
$P(|X_1| \notin \{ 0 \} \cup \{ \psi(K_m):m\in\Z^+ \})=0$.

Let $\e=1$.  We have
$$
        E[\phi(|X_1|)]
         =\sum_{m=1}^\oo 2^{-m} K_m^{-1} \cdot K_m = 1.
$$
It is easy to see that in general $E[\phi(|X_1|)]<\oo$ if and only if (b) holds,
and hence indeed (b) is satisfied in the present case.  So is (a), since $X_1$
is symmetric and so we can put $\mu_n=0$ for all $n$.
Now, let $M(n)=\max \{ m\in \Z^+_0 : K_m< n \}$.  Then, with $T_{1,n}$
as in condition~(iii) of Theorem~\ref{th:main}, and as $K_m\ge 2$ for
all $m\ge 2$
\[
\begin{split}
   T_{1,n}&=\sum_{m=1}^{M(n)} (2^{-m}/K_m) (\psi(K_m))^2 \\
   &= \sum_{m=1}^{M(n)} (2^{-m}/K_m) (K_m \log K_m) \\
   &= \sum_{m=1}^{M(n)} 2^{-m}\log K_m \ge 2^{-M(n)}\log K_{M(n)}.
\end{split}
\]
Note that if $K_m <n \le K_{m+1}$, then $M(n)=m$ and so $T_{1,n}\ge
2^{-m}\log K_m$.  Thus, for $m\ge 1$ we have
\begin{equation}\label{eq:oneterm}
        \sum_{n=K_m+1}^{K_{m+1}}
        n^{-1-1/T_{1,n}} 
        \ge\sum_{n=K_m+1}^{K_{m+1}}
         n^{-1-2^m/\log K_m}.
\end{equation}
Now, for any $K\ge 2$ and $m\in\Z^+$, let $L(K,m)$ be an integer greater than $K$ and
sufficiently large that:
\begin{equation}
        \sum_{n=K+1}^{L(K,m)}
          n^{-1-2^m/\log K}
          \ge \frac12 \sum_{n=K+1}^\oo
           n^{-1-2^m/\log K}.
\end{equation}
Such an $L(K,m)$ exists because the sum on the right hand side
converges.  Observe furthermore that
\begin{equation}\label{eq:threeterm}
\begin{split}
        \sum_{n=K+1}^\oo
           n^{-1-2^m/\log K} \ge C \cdot 2^{-m} (\log K) (K+1)^{-2^m/\log K},
\end{split}
\end{equation}
for an absolute constant $C>0$ independent of $K\ge 2$ and $m\ge 0$. Combining
\eqref{eq:oneterm}--\eqref{eq:threeterm} we see that if $K_{m+1}\ge
L(K_m,m)$, then we have
\begin{equation}\label{eq:the-ineq00}
\begin{split}
        \sum_{n=K_m+1}^{K_{m+1}}
        n^{-1-1/T_{1,n}} 
         &\ge \frac C2 \cdot 2^{-m} (\log K_m) (K_m+1)^{-2^m/\log K_m}
         \\
         &= 2^{-m-1} C (\log K_m) \exp\left(-2^m
\frac{\log(K_m+1)}{\log K_m}\right).
\end{split}
\end{equation}
Inductively choosing the $K_m$ so that for all $m$
we have both
$$
2^{-m-1} (\log K_m) \cdot \exp\left(-2^m
\frac{\log(K_m+1)}{\log K_m}\right) \ge 1,
$$
and $K_{m+1}\ge
L(K_m,m)$, we then find by \eqref{eq:the-ineq00} that
that
\[
        \sum_{n=2}^\oo
         n^{-1-1/T_{1,n}}
        \ge \sum_{m=1}^\oo \sum_{n=K_m+1}^{K_{m+1}}
        n^{-1-1/T_{1,n}}
        \ge C \sum_{m=1}^\oo 1 = \oo.
\]
Hence Corollary~\ref{cor:Sp}'s condition~(c) fails, and so we do have our
desired counterexample satisfying
(a), (b) but not (c), and hence by that Corollary, with inequality
\eqref{eq:sp-conv} also failing.
\end{ex}

Although the answer to Sp\u ataru's question is negative, we do have the 
result under a slightly stronger moment condition than Corollary~\ref{cor:Sp}'s
condition~(b).

\begin{cor}\label{cor:Sp-weak}
Let $X$ be a random variable and let $K$ be any finite constant.
Let $\{ X_{kn} \}_{1\le k\le n}$ be a triangular array of
random variables such that $X_{1n},\dots,X_{nn}$ are $K$-weakly mean dominated
by $X$.
Assume that
\begin{enumerate}
\item[\upn{(a)}] $E[X_{n1}+\dots+X_{nn}]=0$ for all $n$, and
\item[\upn{(b)}] $E\left[\frac{X^2 (\log^+\log^+ |X|)^{1+\delta}}{\log^+ |X|}\right]<\oo$
for some $\delta>0$.
\end{enumerate}
Then,
$$
        \sum_{n=2}^\oo \frac1n P(|S_n| \ge \e (n\log n)^{1/2}) < \oo,
        \quad \forall \e>0.
$$
\end{cor}

\begin{proof} Using Corollary~\ref{cor:wmd} and the same methods as in
Corollary~\ref{cor:Sp} it suffices to show that if condition~(b)
of Corollary~\ref{cor:Sp-weak} holds, then condition~(c) of
Corollary~\ref{cor:Sp} is satisfied with $X$ in place of $X_1$.  
To show this, without loss of generality (rescaling $X$ if necessary)
assume $\e=1$.
Note that if $T_{1,n}=E[X^2\cdot 1_{\{ |X|<a_n
\}}]$, where $a_n=(n\log n)^{1/2}$, then
\begin{equation}\label{eq:star2}
\begin{split}
        T_{1,n} &\le E\left[\frac{X^2 (\log^+\log^+ |X|)^{1+\delta}}{\log^+
               |X|} \right]
                \cdot \frac{\log^+ a_n}{(\log^+\log^+ a_n)^{1+\delta}} \\
                &\le C \frac{\log^+ a_n}{(\log^+\log^+ a_n)^{1+\delta}} \\
                &\le C' \frac{\log^+ n}{(\log^+ \log^+ n)^{1+\delta}},
\end{split}
\end{equation}
where $\delta$ is as in (b), while $C$ and $C'$ are 
strictly positive finite constants independent of $n\ge 2$ (but
dependent on the value of the expectation in (b)).  Hence,
if $N\ge 2$ is sufficiently large that $(\log^+\log^+ n)^{1+\delta}\ge 2C' 
\log\log n$ for all $n\ge N$, then we have:
\[
\begin{split}
        \sum_{n=N}^\oo n^{-1-1/T_{1,n}}
         &=\sum_{n=N}^\oo
            n^{-1} \exp(-(\log n)/T_{1,n}) \\
         &\le \sum_{n=N}^\oo
            n^{-1} \exp(-(\log^+ \log^+
            n)^{1+\delta}/C') \\
         &\le \sum_{n=N}^\oo n^{-1}
          \exp(-2\log\log n) 
         = \sum_{n=N}^\oo \frac1{n (\log n)^{2}} < \oo,
\end{split}
\]
by \eqref{eq:star2} and the choice of $N$.
Thus condition~(c) of Corollary~\ref{cor:Sp} is indeed satisfied 
with $X$ in place of $X_1$.
\end{proof}

\section{Proofs and auxiliary results}\label{sec:proofs}
\begin{proof}[Proof of Proposition~\upref{prp:condA}]
If $\liminf \tau_n>0$ and $\sum_{n=1}^\oo \tau_n
\min(1,nc_n)<\oo$, then there are only finitely many $n\in\Z^+$
for which $nc_n\ge 1$, and hence $\sum_{n=1}^\oo \tau_n nc_n$ must also
converge since it differs from $\sum_{n=1}^\oo \tau_n \min(1,nc_n)$ only
in finitely many terms.

It remains to show that if (a) and (b) hold, then Condition~A.
To do this, suppose $c_n$ is a decreasing sequence such that
\begin{equation}\label{eq:c0}
\sum_{n=1}^\oo \tau_n\min(1,nc_n)<\oo.
\end{equation}
Then, using (b) we have:
\begin{equation}\label{eq:star00}
\begin{split}
        \oo& >\sum_{n=1}^\oo \tau_n\min(1,nc_n) \\
           &\ge \sum_{k=0}^\oo C^{-1}2^k\tau_{2^{k+1}}\min(1,2^k
           c_{2^{k+1}}) \\
           &\ge \frac{1}{2C}\sum_{k=0}^\oo
           2^{k}\tau_{2^{k+1}}\min(1,2^{k+1} c_{2^{k+1}}).
\end{split}
\end{equation}
Now, $\liminf 2^k \tau_{2^{k+1}}>0$ by (a), and hence it follows that only
finitely many of the $\min(1,2^{k+1} c_{2^{k+1}})$ can equal $1$
(since otherwise the right hand side of \eqref{eq:star00} would be
infinite), so
that except for at most finitely many values of $k$, we have $\min(1,2^{k+1}
c_{2^{k+1}})=2^{k+1}c_{2^{k+1}}$.  It thus
follows from \eqref{eq:star00} and (b)
that
\[
\begin{split}
        \oo & > \frac12 \sum_{k=0}^\oo 2^{k+1} \tau_{2^{k+1}} \cdot
        2^{k+1}c_{2^{k+1}} \\
            &\ge \frac{1}{2C}\sum_{k=0}^\oo \sum_{j=0}^{2^{k+1}-1}
              \tau_{j+2^{k+1}} {\frac{j+2^{k+1}}{2}}\cdot c_{j+2^{k+1}} \\
            &= \frac{1}{4C}\sum_{n=2}^\oo \tau_n n c_n,
\end{split}
\]
and the proof is complete.
\end{proof}

\begin{lem}\label{lem:1}  In the setting of Theorem~\upref{th:main}, if the series
$\{ \tau_n \}$ satisfies Condition~A, then condition \eqref{eq:conv}
entails \upn{(i)} and \upn{(ii)}.
\end{lem}

\begin{proof}
Assume \eqref{eq:conv} holds.
Fix $\e>0$ and any
sequence of medians $\mu_n$ of the $S_n$.  Let $M$ be as in
condition~(i).  Then, it is easy to see that for $n\in M$ we have
$P(|S_n|\ge \e a_n)\ge 1/2$.  The convergence of $\sum_{n\in M} \tau_n$
follows immediately from this and \eqref{eq:conv}, and so condition~(i)
holds.

On the other hand, by \eqref{eq:conv} and the remark following the main theorem
in~\cite{Pruss:Klesov}, we have
$$
        \sum_{n=1}^\oo \tau_n \min (1,nP(|X_1|\ge \e a_n)) < \oo
$$
for all $\e>0$.  Condition~(ii) follows immediately from this together
with the fact that
$\{ \tau_n \}$ satisfies Condition~A while $\{ a_n \}$ is increasing.
\end{proof}

\begin{lem}\label{lem:elementary}
Suppose that $\{ \tau_n \}$ and $\{ \rho_n \}$ are positive  and that
$\{ b_n \}$ is increasing and strictly positive.  Let $X$ be any
random variable. If there is a constant $C\in (0,\oo)$
such that for all $n\ge 2$:
\begin{equation}\label{eq:aux-cond-gen}
        b_n^t \sum_{k=n}^\oo \rho_k  \le C \sum_{k=1}^{n-1} k\tau_k,
\end{equation}
then
$$
        \sum_{n=1}^\oo \rho_n E[|X|^t\cdot 1_{\{ |X|<b_n \}}].
\le C        \sum_{n=2}^\oo n\tau_n P(|X|\ge b_n)
$$
\end{lem}

\begin{proof}[Proof of Lemma~\upref{lem:elementary}]
Let $b_0=0$.
Then, by Fubini's theorem and \eqref{eq:aux-cond-gen}:
\[
\begin{split}
  \sum_{n=2}^\oo \rho_n E[|X|^t\cdot 1_{\{ |X|<b_n \}}]
   &\le \sum_{n=2}^\oo \rho_n \sum_{k=1}^{n} b_k^t P(b_{k-1}\le |X|<b_k) \\
   &=\sum_{k=2}^\oo P(b_{k-1}\le |X|<b_k) b_k^t \sum_{n=k}^\oo \rho_n \\
   &\le C\sum_{k=2}^\oo P(b_{k-1}\le |X|<b_k) \sum_{n=1}^{k-1} n\tau_n \\
   &=C \sum_{n=1}^\oo n\tau_n \sum_{k=n+1}^\oo P(b_{k-1}\le |X|<b_k) \\
   &=C \sum_{n=1}^\oo n\tau_n P(|X|\ge b_n).
\end{split}
\]
\end{proof}

\begin{lem}\label{lem:comp} Suppose that $\{ \alpha_n \}$ and $\{ \beta_n \}$ are sequences
in $[0,1]$, that $\{ \tau_n \}$ is a positive sequence, and
that $r\in\Z^+$ is such that\/ $\sum_{n=1}^\oo \tau_n
|\alpha_n-\beta_n|^r<\oo$.  If\/ $\sum_{n=1}^\oo \tau_n \beta_n<\oo$, then
$\sum_{n=1}^\oo \tau_n \alpha_n^r$.
\end{lem}

\begin{proof}
There is a polynomial
$p_r(x,y)$ of degree $r$ with coefficients depending only on $r$ such that
$(x-y)^r=x^r-p_r(x,y)$.  Let $c_r$ be the maximum of $p_r$ over
$[0,1]^2$.  Then:
$$
   \oo > \sum_{n=1}^\oo \tau_n |\alpha_n-\beta_n|^r
    =\sum_{n=1}^\oo \tau_n |\alpha_n^r-\beta_n p_r(\alpha_n,\beta_n)|
    \ge \sum_{n=1}^\oo \tau_n (\alpha_n^r - c_r \beta_n).
$$
If $\sum_{n=1}^\oo \tau_n \beta_n$ converges, then it follows that
$\sum_{n=1}^\oo \tau_n \alpha_n^r$ also converges.
\end{proof}

\begin{lem}\label{lem:the-lemma}
  Let $\{ \tau_n \}$ be a positive sequence and let 
  $\{ b_n \}$ be a strictly positive increasing sequence for $n\ge 1$.
  Let $X$ be a random variable
  such that
\begin{equation}\label{eq:crit}
  \sum_{n=1}^\oo \tau_n n P(|X|\ge b_n) < \oo.
\end{equation}
  Fix $\nu\in[0,\oo)$.  Put $T_k=\sum_{n=1}^k n\tau_n$.
  Suppose that there is a constant $C\in (0,\oo)$ and a $\theta \in
  [1,\oo)$ such that:
\begin{equation}\label{eq:aux-gen-1}
  \frac{b_n^{\nu\theta}}{n^{\theta-1}}\sum_{k=n}^\oo \frac{k^\theta \tau_k}{b_k^{\nu\theta}}
    \le C T_{n-1},
\end{equation}
for all $n\ge 2$, and
\begin{equation}\label{eq:aux-gen-2}
   \frac{k b_n^\nu}{b_k^\nu} \le C T_{n-1},
\end{equation}
whenever $k\ge n\ge 2$.  Then:
$$
\sum_{n=1}^\oo \tau_n \left(\frac{nE[|X|^\nu\cdot 1_{\{ |X_1|<b_n \}}]}
               {b_n^\nu}\right)^{\!\!\theta}<\oo
$$
\end{lem}

The proof is based on methods of Klesov~\cite[Proof of Theorem~4]{Klesov}.

\begin{proof}
Set $X^{(n)}=X\cdot 1_{\{ |X|<b_n \}}$.
Let $b_0=0$ and put $t_n=E[|X^{(n)}|^\nu]$ for $n\ge 0$.  Note that $t_0=0$.
Put $\delta_n=t_n^\theta-t_{n-1}^\theta$ for $n\ge 1$.
For convenience, let $T_1=1$;  redefining $C$ if necessary, we may assume
that \eqref{eq:aux-gen-1} and \eqref{eq:aux-gen-2} also hold for $n=1$.
Then, since $t_n^\theta=\sum_{k=1}^n
\delta_k$, and using Fubini's theorem and \eqref{eq:aux-cond}:
\begin{equation}\label{eq:Klesov-1}
\begin{split}
\sum_{n=1}^\oo \tau_n
        \left(\frac{nE[|X^{(n)}|^\nu]}{b_n^\nu}\right)^{\!\!\theta}
      &= \sum_{n=1}^\oo \tau_n (nb_n^{-\nu})^\theta \sum_{k=1}^n \delta_k \\
      &= \sum_{k=1}^\oo \delta_k \sum_{n=k}^\oo \tau_n (nb_n^{-\nu})^\theta \\
      &\le C\sum_{k=1}^\oo \delta_k k^{\theta-1} b_k^{-\nu\theta} T_{k-1}.
\end{split}
\end{equation}
Let $P_k=P(b_{k-1}\le |X| < b_k)$.  Then, note that
\begin{equation}\label{eq:delta-k}
        \delta_k \le c (t_k-t_{k-1})t_k^{\theta-1} \le cP_k b_k^\nu t_k^{\theta-1},
\end{equation}
where $c$ is a finite constant depending only on $\theta$ and such 
that $x^\theta-y^\theta\le c (x-y)x^{\theta-1}$ whenever $y\le x$.
Now, fix $l\in\Z^+$.  Let $\rho_n=l/b_l^\nu$ for $n=l$ and put $\rho_n=0$ for all other
$n$.  Observe that by \eqref{eq:aux-gen-2} we have
$$
        b_l^\nu \sum_{j=n}^\oo \rho_j  \le C T_{n-1},
$$
for all $n$.
By Lemma~\ref{lem:elementary} (with $t=\nu$) and using \eqref{eq:aux-gen-1} 
together with the definition of the $\rho_n$,
it follows that if $l\ge 2$,
then
\[
\begin{split}
        \frac{lt_l}{b_l^\nu}
        &=\sum_{n=1}^\oo \rho_n t_n
        =\sum_{n=1}^\oo \rho_n E[|X^{(n)}|^\nu] \\
        &\le C \sum_{n=1}^\oo n\tau_n
        P(|X|\ge b_n)
        <\oo,
\end{split}
\]
where the finiteness of the right hand side followed from \eqref{eq:crit}.
Thus, $K\eqdef\sup_{l\ge 1} lt_l/b_l^\nu<\oo$.
Then, $t_k \le K
b_k^\nu/k$, so that \eqref{eq:delta-k} yields:
$$
        \delta_k \le cK^\theta P_k b_k^\nu (b_k^\nu/k)^{\theta-1}=cK \frac{P_k b_k^{\nu\theta}}{k^{\theta-1}}.
$$
Putting this into \eqref{eq:Klesov-1}, and recalling that $T_1=1$, we see that:
\[
\begin{split}
        \sum_{n=1}^\oo \tau_n
        \left(\frac{nE[|X^{(n)}|^\nu]}{b_n^\nu}\right)^{\!\!\theta}
      &\le CcK^\theta\sum_{k=1}^\oo P_k T_{k-1} \\
      &=CcK^\theta\left(P_1+\sum_{k=2}^\oo P_k \sum_{n=1}^{k-1} n\tau_n)\right) \\
      &=CcK^\theta\left(P_1+\sum_{n=1}^\oo n\tau_n \sum_{k=n+1}^\oo P_k\right) \\
      &=CcK^\theta\left(P_1+\sum_{n=1}^\oo n\tau_n P(|X_1|\ge b_n)\right) < \oo,
\end{split}
\]
by Fubini's theorem and \eqref{eq:crit}.
\end{proof}

The following version of the Hoffman-J\o rgensen inequality~\cite{HJ} will be
needed for the proof of Theorem~\ref{th:main} in the case $\theta>1$ and
follows immediately from \cite[Lemma~2.2]{LRJW}.

\begin{lem}\label{lem:HJ} Let $X_1,\dots,X_n$ be independent symmetric random
variables, and let $S_n=X_1+\dots+X_n$.  Then for each $r \in\Z^+$ there
exist finite constants $C_r$ and $D_r$ such
that for all $\lambda\ge 0$ we have:
$$
        P(|S_n|\ge \lambda)
         \le C_r \sum_{k=1}^n P(|X_1|\ge \lambda/(2r))
         + D_r [P(|S_n|\ge \lambda/(2r))]^r.
$$
\end{lem}

Let $\Phi$ be the distribution function of a (0,1) normal random variable.
Recall that a random variable $X$ is {\em symmetric} if and only if $X$
and $-X$ have the same distribution.

\begin{lem}\label{lem:2}
Under the global conditions of Theorem~\upref{th:main},
suppose that condition \upn{(ii)} is satisfied and that
$X_1$ is symmetric.
Then, the following four conditions are equivalent:
\begin{enumerate}
\item[\upn{(a)}] $\sum_{n=1}^\oo \tau_n P(|S_n|\ge \e a_n)<\oo$ for
all $\e>0$;
\item[\upn{(b)}] for all $\gamma\ge 1$ and $\e>0$ we have
$\sum_{n=1}^\oo \tau_n [1-\Phi(\gamma\e a_n/(nT_{\e,n})^{1/2})]<\oo$,
where $T_{\e,n}$ is as in Theorem~\upref{th:main},
and where $\Phi(t/0)=0$ for all $t>0$;
\item[\upn{(c)}] there is an $s\ge 1$ such that for all $\e>0$ we have
$\sum_{n=1}^\oo \tau_n [1-\Phi(\e a_n/(nT_{\e,n})^{1/2})]^s<\oo$;
\item[\upn{(d)}] condition~\upn{(iii)} of Theorem~\upref{th:main} holds.
\end{enumerate}
\end{lem}

The equivalence of (a) and (b) will be the most difficult part to prove,
and will involve a similar method of proof to that in~\cite{Pruss:Riemann},
using a central limit theorem estimate in the present case due to
Nagaev~\cite{Nagaev}.  In the case where $\theta>1$, we will also need
the Hoffman-J\o rgensen inequality based methods of Klesov~\cite[Proof
of Theorem~4]{Klesov}.

\begin{proof}[Proof of Lemma~\upref{lem:2}]
Note that:
\begin{equation}\label{eq:Phi-ineq}
        1-\Phi(x)=\pi^{-1/2} x^{-1} e^{-x^2} (1-O(x^{-2}))
\end{equation}
as $x\to\oo$.  Using the fact that $T_{\gamma \e,n}\le T_{\e,n}$ for 
$\gamma\ge 1$, we see by \eqref{eq:Phi-ineq} that (d) implies (b).
The implication from (b) to (c) is trivial.  Suppose now that (c) holds.
Let $u_{\delta,n}=a_n/(nT_{\delta,n})^{1/2}$.
Then,
\begin{equation}\label{eq:from-b}
        \sum_{n=1}^\oo \tau_n u_{\delta,n}^{-s} e^{-s
         \delta^2 u_{\delta,n}^2} <
        \oo,
\end{equation}
for all $\delta>0$ by (c) and \eqref{eq:Phi-ineq}.  Fix $\e>0$ and let
$\delta=\e/(2s^{1/2})$.  Note that $u_{\delta,n}\ge u_{\e,n}$ and so
\begin{equation}\label{eq:comp-exp}
e^{-\e^2 u_{\e,n}^2} \le K_\e u_{\delta,n}^{-s} e^{-s \delta^2
u_{\delta,n}^2},
\end{equation}
for all $n$, where
$$
        K_\e=\sup_{x\ge 0} x^s e^{(s \delta^2 - \e^2)x^2}
         =\sup_{x\ge 0} x^s e^{-\e^2x^2/2} < \oo.
$$
Then, (d) follows from \eqref{eq:Phi-ineq}--\eqref{eq:comp-exp}, as
desired.  Hence (c) implies (d), so that we have shown that
$\text{(d)}\Rightarrow\text{(b)}\Rightarrow\text{(c)}\Rightarrow\text{(d)}$.

All we now need to prove is the equivalence of (a) and (b).
To do this, assume we are in the setting
of Theorem~\ref{th:main} and that (ii) holds.
Fix $\e>0$.  Let
$X_{k,n}(\e)=X_{k}\cdot 1_{\{ |X_k|<\e a_n \}}$.  Let
$S_n(\e)=X_{1n}(\e)+\dots+X_{nn}(\e)$.  Set $A_n(\e)=\bigcup_{k=1}^n \{
X_k\ne X_{kn}(\e) \}$.
Observe that $S_n=S_n(\e)$ except possibly on $A_n(\e)$, and that
\begin{equation}\label{eq:A-sum}
        \sum_{n=1}^\oo \tau_n P(A_n(\e)) \le \sum_{n=1}^\oo n\tau_n
        P(|X_1|\ge \e a_n) < \oo,
\end{equation}
by (ii).  It follows from \eqref{eq:A-sum} and from the equality of
$S_n$ and $S_n(\e)$ outside $A_n(\e)$ that (a) holds if and only if
\begin{equation}\label{eq:convp}
        \sum_{n=1}^\oo \tau_n P(|S_n(\e)|\ge \e a_n)<\oo,\qquad\forall
        \e>0.
\end{equation}

We now need only show that \eqref{eq:convp} holds if and only if (b)
holds, and we will be done.  Note that $T_{\e,n}=E[(X_{1n}(\e))^2]$.
Fix $\gamma>0$ to be chosen later as needed.  Then,
since all the $X_{kn}(\e)$ and $S_n(\e)$ have mean zero by symmetry, and since
$X_{1n}(\e),\dots,X_{nn}(\e)$ are identically distributed for a fixed $n$, by
Nagaev's central limit theorem estimate~\cite[Theorem~3]{Nagaev} we have:
\begin{equation}\label{eq:Nagaev}
        \left|P(|S_n(\e)|\ge \gamma\e a_n)-2\left[1-\Phi\left(
\frac{\gamma\e a_n}{nT_{\e,n}^{1/2}}\right)\right]\right|
         \le c n \frac{E[|X_{1n}(\e)|^3]}{(\gamma\e a_n)^3},
\end{equation}
for an absolute constant $c<\oo$.

Observe that if \eqref{eq:aux-cond} holds for some $\theta$, then
it also holds for all greater values.
We now have a quick proof if $\theta=1$.  For then, by \eqref{eq:Nagaev}
(with $\gamma=1$)
we
can see that the equivalence of (b) and \eqref{eq:convp} would follow
as soon as we could show that we have
\begin{equation}\label{eq:final-conv}
        \sum_{n=1}^\oo n \tau_n \frac{E[|X_{1n}(\e)|^3]}{a_n^3}<\oo.
\end{equation}
But \eqref{eq:final-conv} follows from the validity of
\eqref{eq:aux-cond} for $\theta=1$ and from condition~(ii) of
Theorem~\ref{th:main}, by an application of Lemma~\ref{lem:elementary}
with $\rho_n=n\tau_n/a_n^3$, $X\equiv X_1$, $t=2$ and  $b_n=\e a_n$.

Suppose now we are working with $\theta>1$, so that
\eqref{eq:aux-aux-cond} also holds.  Let $r$ be an integer greater than
or equal to $\theta$.  I now claim that:
\begin{equation}\label{eq:the-claim}
        \sum_{n=1}^\oo \tau_n
        \left(\frac{nE[|X_{1n}(\e)|^3]}{a_n^3}\right)^{\!\!r} < \oo, \qquad
        \forall \e>0.
\end{equation}

Suppose for now that this has been shown.  If (a) holds, then as noted
before, \eqref{eq:convp} does likewise.  Letting
$\alpha_n=2[1-\Phi(\e a_n/(nT_{\e,n})^{1/2}))]$ and
$\beta_n=P(|S_n(\e)|\ge \e a_n)$, we see that by \eqref{eq:convp}
together with \eqref{eq:the-claim} and \eqref{eq:Nagaev} (with
$\gamma=1$), we do have the conditions of Lemma~\ref{lem:comp}
satisfied, so that $\sum_{n=1}^\oo \tau_n \alpha_n^r < \oo$, and hence
(c) follows, whence (b) follows by the already proved equivalence of
(b), (c) and (d).

Conversely, suppose (b) holds.  This time letting
$\beta_n=2[1-\Phi(\gamma\e a_n/(nT_{\e,n})^{1/2}))]$ and
$\alpha_n=P(|S_n(\e)|\ge \gamma\e a_n)$, using (b), together with
\eqref{eq:Nagaev}, \eqref{eq:the-claim} and Lemma~\ref{lem:comp}, we see
that $\sum_{n=1}^\oo \tau_n \alpha_n^r<\oo$, i.e.,
\begin{equation}\label{eq:res1}
   \sum_{n=1}^\oo \tau_n [P(|S_n(\e)|\ge \gamma\e a_n)]^r<\oo.
\end{equation}
Let $\gamma=(2r)^{-1}$.  By the Hoffman-J\o rgensen inequality
(Lemma~\ref{lem:HJ}), we have:
\[
\begin{split}
  \sum_{n=1}^\oo &\tau_n P(|S_n(\e)|\ge \e a_n) \\
    &\le C_r \sum_{n=1}^\oo \tau_n P(|X_n(\e)|\ge \gamma\e a_n)
     + D_r \sum_{n=1}^\oo \tau_n [P(|S_n(\e)|\ge \gamma\e a_n)]^r.
\end{split}
\]
It is easy to see that the first sum on the right hand side is no greater than $\sum_{n=1}^\oo
\tau_n P(A_n(\gamma\e))$, which converges by \eqref{eq:A-sum}, and the
second converges by \eqref{eq:res1}.  Hence, \eqref{eq:convp} follows,
and as already shown this implies (a).

Hence, all we need to show is that \eqref{eq:the-claim} holds.
Changing finitely many values of $a_n$ and $\tau_n$, we may assume that
\eqref{eq:aux-cond} holds for all $n\ge 2$, and that \eqref{eq:aux-aux-cond}
gives:
\begin{equation}\label{eq:aux-cond-mod}
        \sup_{k\ge n} \frac{k a_n^3}{a_k^3}
        \le C'\sum_{k=1}^{n-1} n\tau_n,
\end{equation}
for a finite $C'$.  We can now apply Lemma~\ref{lem:the-lemma} with $\nu=3$,
$b_n=\e a_n$, and $X\equiv X_1$, using the assumed condition~(ii) of 
Theorem~\ref{th:main} to guarantee \eqref{eq:crit}, and getting
\eqref{eq:aux-gen-1} and \eqref{eq:aux-gen-2} (with an appropriately
chosen constant) from \eqref{eq:aux-cond} and
\eqref{eq:aux-cond-mod} (which hold for all $n\ge 2$ by
assumption), 
respectively.  The Lemma then immediately yields \eqref{eq:the-claim} and
the proof is complete.
\end{proof}

\begin{proof}[Proof of Theorem~\upref{th:main}]
Let $\alpha_1,\alpha_2,\dotsc$ be a sequence of independent Bernoulli
random variables with $P(\alpha_n=1)=P(\alpha_n=-1)=\half$, and with the
sequence independent of $\{ X_n \}_{n=1}^\oo$. Let $X_n'=\alpha_n X_n$,
and put $S_n'=X_1'+\dots+X_n'$.
Note that the primed variables are symmetric.

Suppose first that \eqref{eq:conv} holds and $\{ \tau_n \}$ satisfies
Condition~A.  By Lemma~\ref{lem:1} we have (i) and (ii) holding.  It
remains to show that (iii) holds.    But, if \eqref{eq:conv} holds, it
likewise holds with $S_n'$ in place of $S_n$, as can be seen by
conditioning on the $\{ \alpha_k \}_{k=1}^\oo$ and
using~\cite[Corollary~5]{MS}. But then by the implication
$\text{(a)}\Rightarrow\text{(d)}$ of Lemma~\ref{lem:2} applied to the
primed variables (which are symmetric) it follows that (iii) holds with
$X_1'$ in place of $X_1$.  But (iii) holding for $X_1'$ is equivalent to
it holding for $X_1$, and so (iii) follows.

Conversely, suppose (i), (ii) and (iii) hold.  Evidently (ii) and (iii)
will also
hold with $X_1'$ in place of $X_1$.  Hence by
the implication $\text{(d)}\Rightarrow\text{(a)}$ of Lemma~\ref{lem:2}
as applied to the primed variables we have
\begin{equation}\label{eq:symm-conv}
        \sum_{n=1}^\oo \tau_n P(|S_n'|\ge \e a_n)<\oo, \qquad
        \forall\e>0.
\end{equation}
Now, for any random variable $Y$, let $Y^s=Y-\tilde Y$ be the
symmetrization of $Y$, where $\tilde Y$ is an independent copy of $Y$.
We shall take symmetrizations in such a way that $X_1^s,X_2^s,\dotsc$
are independent and $S_n^s=X_1^s+\dots+X_n^s$ for all $n$.  Now,
since
\[
\begin{split}
        P(|X_k^s|\ge \lambda)&\le P(|X_k|\ge \lambda/2)+P(|\tilde
        X_k|\ge \lambda/2)\\
        &=2P(|2X_k|\ge \lambda)
        =2P(|2X_k'|\ge \lambda),
\end{split}
\]
for all $\lambda\ge 0$, it follows from \cite[Theorem~1]{Pruss:noniid}
that there is an absolute constant $c>0$ such that
$$
        P(|S_n^s|\ge \lambda)\le c P(|2S_n'|\ge \lambda/c)
$$
for all $\lambda\ge 0$.  By \eqref{eq:symm-conv} it then follows that
$$
        \sum_{n=1}^\oo \tau_n P(|S_n^s|\ge \e a_n)<\oo,
        \qquad\forall\e>0.
$$
By standard symmetrization inequalities~\cite[\S17.1.A]{Loeve} it follows
that
\begin{equation}\label{eq:symm-conv2}
        \sum_{n=1}^\oo \tau_n P(|S_n-\mu_n|\ge \e a_n)<\oo,
        \qquad\forall\e>0,
\end{equation}
where $\mu_n$ is the median of $S_n$ occurring in (i).
Now, for any $\e>0$ we have
\begin{multline}\label{eq:for-ref}
        \sum_{n=1}^\oo \tau_n P(|S_n|\ge \e a_n) \\
         \le \sum_{n=1}^\oo \tau_n P(|S_n-\mu_n|\ge \e a_n/2)
           + \sum_{n=1}^\oo \tau_n \cdot 1_{\{ |\mu_n|>\e a_n/2 \}}.
\end{multline}
By condition~(i), the second sum converges, and by
\eqref{eq:symm-conv2}, so does the first, and hence \eqref{eq:conv}
follows.
\end{proof}

\begin{proof}[Proof of Theorem~\upref{th:aux}]
Assume condition~(ii) of Theorem~\ref{th:main} holds.  Fix $\e>0$.
Changing a finite number of values of $\tau_n$ and $a_n$ 
and using \eqref{eq:aux-cond2} and \eqref{eq:aux-aux-cond2}
will let us assume that if we let $b_n=\e a_n$ and $\nu=2$, then 
conditions \eqref{eq:aux-gen-1} and
\eqref{eq:aux-gen-2} of Lemma~\ref{lem:the-lemma} will be verified
for an appropriate constant.
Applying that Lemma and using condition~(ii) of
Theorem~\ref{th:main} shows that:
\begin{equation}\label{eq:t-aux-1}
        \sum_{n=1}^\oo \tau_n n a_n^{-2} T_{\e,n} < \oo,
\end{equation}
where $T_{\e,n}$ is as in Theorem~\ref{th:main}(iii).
Now, using the elementary inequality $e^{-1/x} \le 4 x^2$ which is valid 
for all $x\ge 0$, we see that
\eqref{eq:t-aux-1} entails that
$$
\sum_{n=1}^\oo \tau_n e^{-\e^2 a_n^2 / (nT_{\e,n})} < \oo,
$$
and so (iii) is true.
\end{proof}

\begin{proof}[Proof of Lemma~\upref{lem:Sp}]
Let $a_n=(n\log n)^{1/2}$ for $n\ge 2$.  Assume the conditions of the
Lemma hold.
Without loss of generality $E[X_1]=0$.
Let $Y_n=X_1\cdot 1_{\{|X_1|<a_n\}}$.
Observe that Condition~(b) of Corollary~\ref{cor:Sp} is equivalent
to the claim that $\sum_{n=2}^\oo P(|X_1|\ge a_n)<\oo$, and note
that $\sum_{k=n}^\oo 1/a_k^2 = O(n/a_n^2)$ as $n\to\oo$.
Therefore, by \cite[Theorem~VI.15, p.~225]{Petrov}, we will have
$S_n/a_n\to 0$ almost surely if and only if we
can show that
\begin{equation}\label{eq:show-3}
        \lim_{n\to\oo} \frac{n E[Y_n]}{a_n} \to 0.
\end{equation}
To prove \eqref{eq:show-3}, first let $U_n=X_1-Y_n$.  Observe
that $E[U_n]=-E[Y_n]$ by condition~(a) of Corollary~\ref{cor:Sp}.  Hence, 
we only need to show that $a_n^{-1} n
E[U_n]\to 0$.  Now, since $|U_n|\ge a_n$ whenever $U_n\ne 0$
and as $x^{-1}\log^+ x$ is decreasing in $x>0$, we have
(using the convention that $(0^2/\log^+ |0|)\cdot (\log^+ |0|)/|0|=0$):
\[
\begin{split}
        a_n^{-1} n E[|U_n|]
          &= a_n^{-1} n E[(U_n^2 / \log^+ |U_n|) \cdot (\log^+ |U_n|) /
          |U_n|]
          \\
          &\le a_n^{-1} n E[(U_n^2/\log^+ |U_n|)\cdot (\log^+ a_n)/a_n] \\
          &= (\log^+ a_n) n a_n^{-2} E[(U_n^2/\log^+ U_n)] \to 0,
\end{split}
\]
since $(\log^+ a_n)na_n^{-2}$ is bounded while $E[(U_n^2/\log^+ U_n)]\to
0$ by dominated convergence, as $U_n\to 0$ almost surely and
$U_n^2/\log^+ U_n \le X_1^2/\log^+ X_1$, whereas $X_1^2/\log^+ X_1 \in L^1$ by (b).
\end{proof}

\section*{Acknowledgments}
The author is most grateful for a number of interesting e-mail
conversations on these topics with Professors Aurel Sp\u ataru and
Stephen Montgomery-Smith.  In particular, the author is very grateful to
Professor Montgomery-Smith for permission to use his
Example~\ref{ex:MS}, and to Professor Sp\u ataru for a number of useful
comments on the early drafts of this paper, including pointing out a gap
in the original proof of Corollary~\ref{cor:Sp} and an error in the
original Example~\ref{ex:MS}.

\providecommand{\bysame}{\leavevmode\hbox to3em{\hrulefill}\thinspace}

\end{document}